\documentclass[11pt,psamsfonts,twoside]{article}
\usepackage{amsmath,amssymb,latexsym}
\usepackage{diagrams}
\usepackage{rotating}
\usepackage[mathscr]{eucal}

\newcommand\qed{{\hspace*{\fill}Q.E.D.\vskip12pt plus 1pt}}

\newcommand\sF{{\mathscr F}}

\newcommand\sO{{\mathscr O}}

\def\ca{\operatorname{ca}}

\def\amp{\operatorname{amp}}

\def\lim{\operatorname{lim}}

\def\Proj{\operatorname{Proj}}

\def\codim{\operatorname{codim}}

\def\cd{\operatorname{cd}}
\def\Coh{\operatorname{Coh}}

\newcommand\pn[1]{{\mathbb P}^{#1}}

\newcommand\proof{\noindent{\em Proof.}\ \ }

\newtheorem{theorem}{Theorem}[section]
\newtheorem{lemma}[theorem]{Lemma}
\newtheorem{corollary}[theorem]{Corollary}

\newtheorem{question}[theorem]{Question}

\newtheorem{definition}[theorem]{Definition}
\newtheorem{rem}[theorem]{Remark}

\newtheorem{pargrph}[theorem]{}
\newtheorem{examp}[theorem]{Example}

\newtheorem{MM}[theorem]{ }
\newtheorem{res}[theorem]{Remarks}

\makeatletter
\ifnum\@ptsize=0  \fi
\ifnum\@ptsize=2  \fi
\renewcommand{\qed}{\hfill $\square$}

\newenvironment{rem*}{\begin{rem}\em}{\end{rem}}
\newenvironment{rems*}{\begin{res}\em}{\end{res}}
\newenvironment{example*}{\begin{examp}\em}{\end{examp}}
\newenvironment{definition*}{\begin{definition}\em}{\end{definition}}
\newenvironment{question*}{\begin{question}\em}{\end{question}}
\newenvironment{MM*}{\begin{MM}\em}{\end{MM}}

\newenvironment{prgrph*}[1]{\indent\begin{pargrph}{\bf #1.}\em\
}{\end{pargrph}}

\begin{document}

\title{A connectedness theorem for products of \\ weighted projective spaces}
\author{Lucian B\u adescu and Flavia Repetto}
\date{}

\maketitle

\begin{abstract} We prove a connectedness result for products of weighted projective spaces.\end{abstract}

\section*{Introduction}
Let $P$ be a projective irreducible variety and let $f\colon X\to P\times P$ be a morphism from a complete irreducible variety $X$ over an algebraically closed field $k$. Denote by $\Delta$ the diagonal of $P\times P$. Then one can ask under which conditions the inverse image $f^{-1}(\Delta)$ is connected (resp. non-empty). Here by a connected scheme we shall mean a non-empty scheme whose underlying topological space is connected. The first result in this direction is the famous theorem of Fulton-Hansen which states that the answer to this question is affirmative if 
$P=\mathbb P^n$ and if  $\dim(f(X))>n$ (resp.  if  $\dim (f(X))\geq n$) (see 
\cite{FH}, cf. also \cite{De} or \cite{FL}). This result has a lot of interesting geometric consequences (loc. cit.).

The theorem of Fulton-Hansen has been subsequently generalized in various directions. First Hansen proved that
$f^{-1}(\Delta)$ is connected (resp. non-empty) if $P=\mathbb{G}(d,n)$ is the Grassmann variety of linear $d$-spaces in 
$\pn n$ and if $\codim_{P\times P}(f(X))<n$ (resp. $\codim_{P\times P}(f(X))\leq n$) see \cite{Han}. 

A quite general connectivity result was proved by Faltings (see \cite{F}) if $P$ is an arbitrary projective rational homogeneous space over an algebraically closed field of characteristic zero. In a slightly improved formulation (via 
\cite{Go}) it states that $f^{-1}(\Delta)$ is connected (respectively, non-empty) if $\codim_{P\times P}(f(X))<\ca(P)$ (resp. 
$\codim_{P\times P}(f(X))\leq \ca(P)$) (see \cite{F} and \cite{Go}). Here $\ca(P)$ is the coampleness of $P$, defined as follows (see \cite{Go}). 
Since $P$ is a homogeneous space, the tangent bundle $T_P$ of $P$ is generated by its global
sections; this implies that the tautological line bundle ${\mathcal O}_{{\mathbb P}(T_P)}(1)$ is also generated by its global sections. Then one defines the ampleness, $\amp(P)$, of $P$ as the maximum fiber dimension of the morphism 
$\varphi\colon{\mathbb P}(T_P)\to {\mathbb P}^N$ associated to the complete linear system $|{\mathcal
O}_{{\mathbb P}(T_P)}(1)|$. Then one defines the coampleness $\ca(P)$ of $P$ by $\ca(P):=\dim(P)-\amp(P)$. A result of Goldstein (\cite{Go}) asserts that $\ca(P)\geq r$, where $r$ is the minimum of ranks of the simple factors of the linear algebraic group $G$ acting transitively on $P$; in particular, $\ca(P)\geq 1$. For example, $\ca(\mathbb P^n)=n$ (or more generally, $\ca( \mathbb{G}(d,n))=n$) and $\ca(\mathbb P^m\times\mathbb P^n)=n$, if $m\geq n\geq 1$. In particular, if $P=\mathbb P^m\times\mathbb P^n$, then Faltings' result asserts that $f^{-1}(\Delta)$ is connected (respectively, non-empty) if $\codim_{P\times P}(f(X))<n$ (resp. $\codim_{P\times P}(f(X))\leq n$).

More recently Debarre  proved two interesting connectivity results, when $P$ is either a product of projective spaces or a Grassmann variety (see \cite{D}). By imposing some extra (geometric) conditions his results go beyond
the bound $\codim_{P\times P}(f(X))<n$, in the case $P=\mathbb P^m\times\mathbb P^n$, or below Hansen's bound  in the case $P=\mathbb{G}(d,n)$. 

Finally, in \cite{B1} (see also \cite{B}, Theorem 7.14) the first named author showed (among other things)  that Fulton-Hansen connectivity result still holds if $P=\mathbb P^n(e)$ is an $n$-dimensional weighted projective space of weights $e=(e_0,e_1,\ldots,e_n)$ ($e_i\geq 1$, $i=0,1,\ldots,n$), provided that 
$\dim(f(X))>n$. 

To state our main result, let $P$ denote the product of the weighted projective spaces $\pn m(e)$ and $\pn n(h)$ of weights 
$e=(e_0,e_1,\ldots, e_m)$ and $h=(h_0,h_1,\ldots,h_n)$, with $e_i,h_j\geq 1$,  $i=0,\ldots,m$ and $j=0,\ldots,n$. Let $f\colon X\to P\times P$ be a morphism from a complete irreducible variety $X$. Denote by $X_{13}\subseteq\mathbb P^m(e)\times\mathbb P^m(e)$ (resp.  by $X_{24}\subseteq\mathbb P^n(h)\times\mathbb P^n(h)$) the image of $f(X)$ under the projection 
$p_{13}$  of $P\times P=\mathbb P^m(e)\times\mathbb P^m(e)\times\mathbb P^m(e)\times\mathbb P^n(e)$ onto $\mathbb P^m(e)\times\mathbb P^n(e)$
(resp. under the projection $p_{24}$ onto $\mathbb P^n(h)\times\mathbb P^n(h)$).
The aim of this paper is to prove the following:

\medskip

{\bf Theorem.} (Theorem \ref{CT} below) {\em In the above notation let $P=\pn m(e)\times\pn n(h)$ be the product of the weighted  projective spaces $\pn m(e)$ and $\pn n(h)$ over an algebraically closed field $k$, let $\Delta$ be the diagonal of $P\times P$ and set $a:=\max\{m+\dim(X_{24}),n+\dim(X_{13})\}$. If  $\dim(f(X))\geq a$ then $f^{-1}(\Delta)$ is non-empty, and if  $\dim(f(X))> a$ then $f^{-1}(\Delta)$ is connected.}

\medskip

In the case when $P=\pn m\times\pn n$ is a product of ordinary projective spaces in Debarre's result (see \cite{D}, Th\'eor\`eme 2.2 and Lemme 2.1) the hypotheses are  the following:
\begin{enumerate}
\item[i)] $\dim(X_{13})>m$ and $\dim(X_{24})>n$ (resp. $\dim(X_{13})\geq m$ and $\dim(X_{24})\geq n$), and
\item[ii)]  $\dim(f(X))>m+n$ (resp. $\dim(f(X))\geq m+n$).
\end{enumerate}
As  is easily seen the hypothesis in our theorem implies i) and ii), and therefore our result,  in the case of a product of ordinary projective spaces, is weaker than Debarre's result.

The above theorem has the following consequence:

\medskip

{\bf Corollary.} {\em  Let $f:X\to P\times P$ be a morphism, with $X$ a complete irreducible variety, $P=\pn m(e)\times \pn n(h)$ and $m\geq n\geq 1$, such that $\codim_{P\times P}(f(X))< n$ $($resp. $\codim_{P\times P}(f(X))\leq n)$. Then $f^{-1}(\Delta)$ is connected $($resp. non-empty$)$.}

\medskip

Notice that this corollary recovers the above mentioned result of Faltings  in the case  when the homogeneous space $P$ is a product of ordinary projective spaces over a field of characteristic zero.

One motivation for proving connectivity theorems  for products of weighted projective spaces lies in the fact that there are several interesting varieties which can be embedded in a nicer way in some weighted projective spaces than in any ordinary projective space. Our approach is completely different from Debarre's (see \cite{D}). In fact, Debarre proceeds as in \cite{FH} 
(see also \cite{FL}) proving first a Bertini theorem for general linear sections and then deduces the connectivity of all linear sections by an argument of degeneration. This method does not seem to apply in the case of weighted projective spaces. Our basic technical ingredient is Hartshorne's cohomological dimension (see \cite{Ha3},  cf. also \cite{Ha2}).

The paper is organized as follows.
In the first section, after recalling the definition  and a few basic facts about the cohomological dimension  of an algebraic variety and about  weighted projective spaces, we present some preparatory material that will be needed in section 2. In the second section we prove the theorem stated above. 

Unless otherwise specified, we shall use the standard terminology and notation in algebraic geometry.
We shall work over an algebraically closed ground field $k$ of arbitrary characteristic.

\section{Background material}\label{first}\addtocounter{subsection}{1}

\begin{definition*}
Let $Z$ be an irreducible algebraic variety over $k$. We shall denote by $\Coh (Z)$ the category of all coherent sheaves of $\sO _Z$-modules. According to Hartshorne \cite{Ha3} (see  also \cite{Ha2}), we define the {\it{cohomological dimension of $Z$}}, $\cd (Z)$, by the following
$$\cd (Z)=\max \{n\geq 0\  | \ \exists \;\sF\in \Coh (Z) \;{\rm{ such\;  that}} \;H^n(Z,\sF)\neq 0\}.$$
\end{definition*}
By a general result of Grothendieck one has $0\leq \cd (Z) \leq \dim (Z)$.
Moreover, by  Serre's criterion of affiness,  $\cd (Z)=0$ if and only if $Z$ is affine.

Now we recall some basic results involving cohomological dimension which are going to be used in the next section. An important result due to Hartshorne and Lichtenbaum (see \cite{Gr65}, \cite{Ha3}, or also  \cite{K} for a more elementary proof) is the following:

\begin{theorem}[Hartshorne-Lichtenbaum]\label{HL}
Let $X$ be an irreducible projective variety of dimension $n\geq 1$, let $Y$ be a closed subset of $X$ and set $U:=X\setminus Y$. Then
 $\cd (U)\leq n-1\  if\  and\  only\ if\ Y\neq\varnothing.$
\end{theorem}

Another  important result is the following (if  $X$ is smooth, see  \cite{Ha2}, Corollary 3.9; for the general case see \cite{B1}, appendix, or also \cite{B}, Theorem 7.6).

\begin{theorem}\label{7.6}
In the notation of Theorem $\ref{HL}$ assume that $X$ is an irreducible projective variety of dimension $n\geq 2$ such that $\cd (U)\leq n-2$. Then $Y$ is connected.\end{theorem} 

The following corollary of Theorems \ref{HL} and \ref{7.6} is going to be our basic technical tool in proving the main result stated in the introduction.

\begin{corollary}\label{universal} Let $Y$ be a closed subset of a projective irreducible variety $X$ of dimension $n\geq 2$ and set $U:=X\setminus Y$.  Let $f\colon X'\to X$ be a finite surjective morphism, with $X'$ irreducible. Then $f^{-1}(Y)$ is non-empty if $\cd(U)\leq n-1$, and connected  if $\cd(U)\leq n-2$.\end{corollary}

\proof Set $U':=f^{-1}(U)$. Then $f|U' \colon U' \to U$ is a finite surjective morphism, whence by a general simple property (see  \cite{Ha3}, Proposition 1.1), $\cd(U')=\cd(U)$. Then the conclusion follows from Theorems \ref{7.6} and \ref{HL} applied to $(X', f^{-1}(Y))$.\qed

\medskip

We shall record for further use the following two well known elementary results:

\begin{lemma}\label{lex} $(${\em \cite{Ha2}, Ex. 3.11 page 103}$)$ For every quasi-projective varieties $V$ and $W$ one has $\cd (V\times W)=\cd (V)+\cd (W)$.\end{lemma}

\begin{lemma}\label{irrid}
Let $\psi:V\to W$ be an open morphism of quasi-projective varieties. If $W$ and all fibers of $\psi$ are irreducible, then $V$ is also irreducible.
\end{lemma}

Finally we shall also need the following result:

\begin{lemma}\label{IV-} $(${\em \cite{IV}, Lemma 4.7}$)$ Let $U_1$ and $U_2$ be two open subsets of a quasi-projective irreducible variety $Z$ such that $Z=U_1\cup U_2$. Assume that  for some integer $a\geq 1$, 
$\cd (U_i)\leq a$  for $i=1,2$, and $\cd (U_1\cap U_2)\leq a-1$. Then $\cd (Z)\leq a$.\end{lemma}

Now we briefly recall the definition and some basic facts on weighted projective spaces.
Let $k[T_0,\ldots ,T_n]$ be the polynomial $k$-algebra  in $n+1$ variables $T_0,\ldots ,T_n$ (with $n\geq 1$). An $(n+1)$-uple $(e_0,\ldots ,e_n)\in\mathbb{Z}^{n+1}$ of positive integers is called a  {\em system of weights}. Given a system of weights $e=(e_0,\ldots ,e_n)$, grade $k[T_0,\ldots ,T_n]$ by the conditions 
$\deg(T_i)=e_i$, for every $i=0,\ldots ,n$. In this way we get a finitely generated graded $k$-algebra (depending of 
$e=(e_0,\ldots ,e_n)$), and set 
$$\pn n (e):= \pn n (e_0,\ldots ,e_n):=\Proj (k[T_0,\ldots ,T_n]).$$ 
Then $\pn n (e)$ is  a projective irreducible, normal variety of dimension $n$, with quotient singularities (hence Cohen-Macaulay), which is called the {\em weighted projective space of weights} $e=(e_0,\ldots ,e_n)$.
A point of $\pn n (e)$ can be given by a collection $[a_0,\ldots,a_n]$ of ordered $n+1$ elements $a_0,\ldots,a_n\in k$, not all of them zero. Another collection $[b_0,\ldots,b_n]$ of $n+1$ elements of $k$ defines the same point of $\pn n(e)$ if and only if there exists a non-zero element $t\in k$ such that $b_i=a_it^{e_i}$, for every $i=0,\ldots,n$. Of course, $\pn n (1,\ldots ,1)$ coincides with the usual projective space $\pn n$. 
For the basic properties of weighted projective spaces see \cite{Do}, or also \cite{BR}.

\begin{definition*} \label{weighted}(See \cite{B1}, or also \cite{B}, page 62).\label{weight1} Let $\pn n (e)$ be the weighted projective space of weights 
$e=(e_0,\ldots ,e_n)$, and let $\pn {2n+1}(e,e)=\Proj(k[T_0,\ldots ,T_n,T'_0,\ldots,T'_n])$ be the weighted projective space of weights $(e,e)=(e_0,\ldots,e_n ;e_0,\ldots ,e_n)$, where $T_0,\ldots,T_n,T'_0,\ldots,T'_n$ are $2n+2$ independent indeterminates over $k$ such that $\deg(T_i)=\deg(T'_i)=e_i$, $i=0,\ldots,n$. Then the canonical inclusions 
$$k[T_0,\ldots ,T_n]\subset k[T_0,\ldots ,T_n,T'_0,\ldots,T'_n]\;\;\text{and}\;\;k[T'_0,\ldots,T'_n]\subset k[T_0,\ldots ,T_n,T'_0,\ldots,T'_n]$$ 
are homomorphisms of graded $k$-algebras which define two rational maps
$$\pi^i_n(e)\colon\pn{2n+1}(e,e)\dashrightarrow\pn n(e),\;\;i=1,2.$$
In fact, $\pi_n^1(e)$ is defined precisely in the complement of $L^n_1(e):=V_+(T_0,\ldots,T_n)$ and $\pi_n^2(e)$ -- in the complement of $L^n_2(e):=V_+(T'_0,\ldots,T'_n)$. Then $\pi_n^1(e)$ and $\pi_n^2(e)$ yield the rational map
$$\pi_n(e)\colon\pn {2n+1}(e,e)\dashrightarrow\pn n(e)\times\pn n(e)$$
which is defined precisely in the open subset $U_n(e):=\pn{2n+1}(e,e)\setminus(L^n_1(e)\cup L^n_2(e))$,
with $L^n_1(e)\cap L^n_2(e)=\varnothing$. In weighted coordinates the rational map $\pi_n(e)$ is defined by
$$\pi_n(e)([t_0,\ldots,t_n,t'_0,\ldots,t'_n])=([t_0,\ldots,t_n],[t'_0,\ldots,t'_n]),$$
$\forall [t_0,\ldots,t_n,t'_0,\ldots,t'_n]\in\mathbb P^{2n+1}(e,e)$.
 We call the morphism $\pi_n(e)\colon U_n(e)\to\pn n(e)\times\pn n(e)$  the {\em weighted join construction}, generalizing the usual join construction used by Deligne in \cite{De} to simplify the proof of Fulton-Hansen connectivity theorem. If $e=(1,\ldots,1)$ we shall simply write $L_i^n$ for $L_i^n(1,\ldots,1)$, $i=1,2$, $\pi_n$ for $\pi_n(1,\ldots,1)$ and $U_n$ for $U_n(1,\dots,1)$. \end{definition*} 

There exists a canonical finite surjective morphism $g:\pn n\to\pn n(e)$.
Indeed, if $\pn n=\Proj (k[X_0,\ldots,X_n])$ with $\deg(X_i)=1$, for $i=0,\ldots,n$, and $\pn n (e)=\Proj(k[T_0,\ldots ,T_n])$ with $\deg(T_j)=e_j$, for $j=0,\ldots,n$, then there is a graded homomorphism of graded rings $\varphi:k[T_0,\ldots ,T_n]\to k[X_0,\ldots,X_n]$,
defined by $\varphi (T_i):=X_i^{e_i}$, for every $i=0,\ldots,n$, which yields the finite surjective morphism $g:\pn n\to\pn n(e).$

We have the following commutative diagram
\begin{diagram}
U_n&\rTo^{g'}&U_n(e)\\
\dTo ^{\pi_n}&  &\dTo _{\pi _n(e)}\\
\pn n\times\pn n&\rTo^{g\times g}&\pn n(e)\times\pn n(e)
\end{diagram}
which is easily seen to be cartesian, and in particular, every fiber of $\pi_n(e)$ is isomorphic to the multiplicative group $\mathbb G_m=k\setminus\{0\}$. As is well known (and easy to see) the morphism $\pi_n$  is a locally trivial $\mathbb{G}_m$-bundle over $\pn n\times\pn n$. In particular,  $\pi_n$ is flat, whence by  \cite{J}, Proposition 2.7, an open morphism.

\begin{lemma}\label{open}
Consider the following cartesian diagram of algebraic varieties:
\begin{diagram}
X'&\rTo^{\alpha}&Y'\\
\dTo ^{\beta}&  &\dTo _{\gamma}\\
X&\rTo^{\delta}&Y
\end{diagram}
with $\delta$ a proper surjective morphism.
If $\beta$ is an open morphism, then $\gamma$ is also open.
\end{lemma}

\proof Since $\delta$ is proper and surjective, so is $\alpha$. Let $V$ be an open subset of $Y'$. We want to prove that $\gamma(V)$ is open in $Y$.  To this extent observe that since the above diagram is cartesian one has
$\beta (\alpha ^{-1}(V))=\delta^{-1}(\gamma(V))$.
Moreover, since $\delta$ is proper and surjective,  the topology on $Y$ is the quotient topology of $X$. Thus $\gamma (V)$ is open in $Y$ if and only if $\delta ^{-1}(\gamma (V))$ is open in $X$.
But  $\delta ^{-1}(\gamma (V))=\beta (\alpha ^{-1}(V))$  is an open subset of $X$ because $V$ is open in $Y'$ and $\beta$ is an open morphism  by hypothesis.
\qed

\begin{corollary}\label{open0} The morphism $\pi_n(e)\colon U_n(e)\to\pn n(e)\times\pn n(e)$ of Definition $\ref{weight1}$ is open. \end{corollary}

\section{The connectedness theorem}\label{second}

Consider the systems of positive weights $e=(e_0,\ldots,e_m)$, $h=(h_0,\ldots,h_n)$, with $m,n\geq 1$,
and the weighted projective spaces $\pn m(e)=\Proj(k[T_0,...,T_m])$ and $\pn n(h)=\Proj(k[U_0,\ldots,U_n])$, with $T_0,...,T_m$ and $U_0,\ldots,U_n$ independent indeterminates of weights
$\deg(T_i)=e_i$, $i=0,...,m$, and $\deg(U_j)=h_j$, $j=0,...,n$. According to the previous section we also consider
$\mathbb P^{2m+1}(e,e)=\Proj(k[T_0,\ldots,T_m,T'_0,\ldots,T'_m])$ and 
$\mathbb P^{2n+1}(h,h)=\Proj(k[U_0,\ldots,U_n,U'_0,\ldots,U'_n])$.

Denote by $X_{13}\subseteq\mathbb P^m(e)\times\mathbb P^m(e)$ (resp. $X_{24}\subseteq\mathbb P^n(h)\times\mathbb P^n(h)$) the image of $X'=f(X)$ under the projection $p_{13}$  of $P\times P=\mathbb P^m(e)\times\mathbb P^n(h)\times\mathbb P^m(e)\times\mathbb P^n(h)$ onto $\mathbb P^m(e)\times\mathbb P^m(e)$
(resp. under the projection $p_{24}$ of $P\times P$ onto $\mathbb P^n(h)\times\mathbb P^n(h)$). Then we prove the following:

\begin{theorem}\label{CT} In the above notation let $P=\pn m(e)\times\pn n(h)$ be the product of the weighted  projective spaces $\pn m(e)$ and $\pn n(h)$ over an algebraically closed field $k$, let $\Delta$ be the diagonal of $P\times P$ and set $a:=\max\{m+\dim(X_{24}),n+\dim(X_{13})\}$. If  $\dim(f(X))\geq a$ then $f^{-1}(\Delta)$ is non-empty, and if  $\dim(f(X))> a$ then $f^{-1}(\Delta)$ is connected.
\end{theorem}

\proof 
Set $X':=f(X)\subseteq P\times P$ and consider the Stein factorization of $f=f_2\circ f_1$, with $f_1:X\to W$ a proper surjective morphism with connected fibers and $f_2:W\to X'$ a finite morphism (see 
\cite{Gr66}, III). Thus we may assume that $f$ is a finite morphism.

According to  the notation of the previous section we set
$$U_m(e):=\pn {2m+1}(e,e)\setminus (L_1^m(e)\cup L_2^m(e))\:\:\text{and}\;\;U_n(h):=\pn {2n+1}(h,h)\setminus (L_1^n(h)\cup L_2^n(h)).$$ 
Since $T_i-T_i'$ (respectively $U_j-U_j'$) is a homogeneous element  of degree $e_i$, $i=0,\ldots ,m$ (respectively $h_j$, $j=0,\ldots ,n$), it makes sense to consider also the closed subschemes 
$$H_m(e):=V_+(T_0-T_0',\ldots ,T_m-T_m') \quad {\rm{and}}\quad  H_n(h):=V_+(U_0-U_0',\ldots ,U_n-U_n')$$ 
 of $\pn {2m+1}(e,e)$ and $\pn {2n+1}(h,h)$ respectively.
Clearly $H_m(e)\subseteq U_m(e)$ and $H_n(h)\subseteq U_n(h)$, whence  
$$H_m(e)\times H_n(h)\subseteq U:=U_m(e)\times U_n(h)\subseteq \pn {2m+1}(e,e)\times \pn {2n+1}(h,h).$$
Consider the rational map
$$\pi\colon\pn {2m+1}(e,e)\times \pn {2n+1}(h,h)\dashrightarrow P\times P=\pn m(e)\times \pn n(h)\times\pn m(e)\times \pn n(h)$$
defined by
$$\pi([t_0,\ldots ,t_m, t_0',\ldots ,t_m'],[u_0,\ldots ,u_n,u_0'\ldots ,u_n'])=$$
$$= ([t_0,\ldots ,t_m],[u_0,\ldots ,u_n],[ t_0',\ldots ,t_m'],[u_0',\ldots ,u_n']),$$
$\forall\ ([t_0,\ldots ,t_m, t_0',\ldots ,t_m'],[u_0,\ldots ,u_n,u'_0\ldots ,u_n'])\in\pn {2m+1}(e,e)\times \pn {2n+1}(h,h)$.
Actually, modulo the canonical isomorphism $$\pn m(e)\times \pn n(h)\times\pn m(e)\times \pn n(h)\cong\pn m(e)\times \pn m(e)\times\pn n(h)\times \pn n(h),$$
$\pi$ is nothing but the product $\pi_m(e)\times\pi_n(h)$ of the rational maps $\pi_m(e)\colon\pn {2m+1}(e,e)\dasharrow \pn m(e)\times \pn m(e)$ and $\pi_n(h)\colon\pn {2n+1}(h,h)\dasharrow \pn n(h)\times \pn n(h)$ of Definition \ref{weight1}.
Then  the  map $\pi$  is defined precisely in the open subset $U=U_m(e)\times U_n(h)$, and all the fibers of $\pi$ are isomorphic to $\mathbb{G}_m\times\mathbb{G}_m$ (Corollary \ref{open0}).
Furthermore, with arguments similar to the proof of Corollary \ref{open0} it is easy to see that the
morphism $\pi\colon U\to P\times P$ is open. It is clear the restriction map $\pi |(H_m(e)\times H_n(h))$ defines an isomorphism $H_m(e)\times H_n(h)\cong \Delta$.

Set $X'=f(X)\subseteq P\times P$ and consider the commutative diagram

\begin{diagram}
U_{X'}\cap(H_m(e)\times H_n(h))&\hookrightarrow&U_{X'}&\hookrightarrow&U\\
\dTo ^{\cong}&  &\dTo ^{\pi _{X'}}& &\dTo _{\pi}\\
X'\cap \Delta&\hookrightarrow&X'&\hookrightarrow&P\times P
\end{diagram}
with $U_{X'}:=\pi^{-1}(X')$ and $\pi_{X'}:=\pi|U_{X'}$ (the restriction of $\pi$ to $U_{X'}$).

Set $U'_{X'}:=U_{X_{13}}\times U_{X_{24}}\subseteq\mathbb P^{2m+1}(e,e)\times\mathbb P^{2n+1}(h,h)$. Since $X'$ is a closed subset of $X_{13}\times X_{24}$,  $U_{X'}$ is a closed subset of 
$U'_{X'}$ of dimension $d:=\dim(X')+2$. Moreover, the theorem on the dimension of fibers shows that the hypothesis that
$\dim(f(X))>a$ (resp. if $\dim(f(X))\geq a$) implies
\begin{equation}\label{i}
\dim(X_{13})>m\;\text{and}\;\dim(X_{24})>n\;\text{(resp.}\;\dim(X_{13})\geq m\;\text{and}\;\dim(X_{24})\geq n\text{)}
\end{equation}

Let $Y_{13}$ be the closure of $U_{X_{13}}$ in $\mathbb P^{2m+1}(e,e)$ and  
$Y_{24}$ the closure of $U_{X_{24}}$ in $\mathbb P^{2n+1}(h,h)$. Then  $Y_{13}\times Y_{24}$ is the closure of 
$U'_{X'}=U_{X_{13}}\times U_{X_{24}}$ in $\mathbb P^{2m+1}(e,e)\times\mathbb P^{2n+1}(h,h)$. Moreover,  
$U'_{X'}\cap(H_m(e)\times H_n(h))=(Y_{13}\times Y_{24})\cap(H_m(e)\times H_n(h))$ and 
$(Y_{13}\times Y_{24})\setminus(H_m(e)\times H_n(h))=V'_1\cup V'_2$, with
$V_1':=(Y_{13}\setminus H_m(e))\times Y_{24}$ and $V'_2:=Y_{13}\times(Y_{24}\setminus H_n(h))$.
Using Lemma \ref{lex} and \eqref{i} we get
$$\cd(V'_1)\leq m+\dim(Y_{24})=m+1+\dim(X_{24}),\;\;\cd(V'_2)\leq n+\dim(Y_{13})=n+1+\dim(X_{13}). $$ 
On the other hand, using our hypotheses, the equality $V'_1\cap V'_2=(Y_{13}\setminus H_m(e))\times(Y_{24}\setminus H_n(h))$ and Lemma \ref{lex}, we get
$$\cd(V'_1\cap V'_2)\leq m+n<\max\{m+1+\dim(X_{24}),\;n+1+\dim(X_{13})\}=a+1.$$
Therefore by Lemma \ref{IV-} we get
$$\cd((Y_{13}\times Y_{24})\setminus(H_m(e)\times H_n(h)))\leq a+1.$$
Let $Y$ denote the closure of $U_{X'}$ in  $\mathbb P^{2m+1}(e,e)\times\mathbb P^{2n+1}(h,h)$. Clearly, 
$\dim(Y)=\dim(U_{X'})=d$. Since $Y\setminus(H_m(e)\times H_n(h))$ is a closed subset of $(Y_{13}\times Y_{24})\setminus(H_m(e)\times H_n(h))$, we get $\cd(Y\setminus(H_m(e)\times H_n(h)))\leq\cd((Y_{13}\times Y_{24})\setminus(H_m(e)\times H_n(h)))$, and thus the last inequality yields
\begin{equation}\label{final}\cd(Y\setminus(H_m(e)\times H_n(h)))\leq a+1.\end{equation}
Therefore  \eqref{final} yields
\begin{equation}\label{final'}\cd(Y\setminus(H_m(e)\times H_n(h)))\leq\begin{cases} d-2,\;\text{if}\;\dim(f(X))>a\\
d-1,\;\text{if}\;\dim(f(X))\geq a
\end{cases}.\end{equation}

On the other hand, consider the cartesian diagram 
\begin{diagram}
\overline{X}=X\times _{X'}U_{X'}&\rTo ^{\overline{f}}&U_{X'}\\
\dTo ^{\pi _{X}}&  &\dTo _{\pi _{X'}}\\
X&\rTo ^f&X'
\end{diagram}

\medskip

{\em Claim} 1. The morphism $\pi_X$ is open.

\medskip

To prove the claim we first observe that it is obvious if $e=(1,\ldots,1)$ and $h=(1,\ldots,1)$, i.e. if
$P$ is a product of ordinary projective spaces. In fact, in this case the morphism $\pi\colon U\to P\times P$ is a locally trivial $\mathbb G_m\times\mathbb G_m$-bundle, and in particular, $\pi$ is a flat morphism. It follows that $\pi_X$ is also flat because flatness is preserved under base change, whence $\pi_X$ is open (see e.g. \cite{J}, Proposition 2.7).
In the general case, let $g(e)\colon\mathbb P^m\to\mathbb P^m(e)$ and $g(h)\colon\mathbb P^n\to\mathbb P^n(h)$ be the canonical finite morphisms, and set 
$P_1:=\mathbb P^m\times\mathbb P^n$, $u:=g(e)\times g(h)\times g(e)\times g(h)\colon P_1\times P_1\to P\times P$, $X'_1:=u^{-1}(X')$ and $u_{X'}:=u|X'_1\colon X'_1\to X'$.
Consider the following commutative diagram 
\begin{diagram}
\overline{X}_1&&\rTo^{\overline{f}_1}&&U_{X'_1}&&\\
&\rdTo_{\pi_{X_1}}&&&\vLine&\rdTo^{\pi_{X'_1}}&\\
\dTo^{u_{\overline{X}}}&&X_1&\rTo^{f_1}&\HonV&&X_1'\\
&&\dTo^{u_X}&&\dTo_{u'}&&\\
\overline{X}&\hLine&\VonH&\rTo^{\overline{f}}&U_{X'}&&\dTo_{u_{X'}}\\
&\rdTo_{\pi_X}&&&&\rdTo^{\pi_{X'}}&\\
&&X&&\rTo^f&&X'\\
\end{diagram}
in which $\overline{X}:=X\times_{X'}U_{X'}$, $X_1:=X\times_{X'}X'_1$, $U_{X'_1}:=U_{X'}\times_{X'}X_1'$, and $\overline{X}_1:=X_1\times_{X'_1}U_{X'_1}$. An easy diagram chase shows that  the left vertical square is cartesian, i.e. $\overline{X}_1=\overline{X}\times_XX_1$. Now, since $P_1$ is a product of ordinary projective spaces, by what we have said above, the morphism $\pi_{X_1}$ is flat and hence open. Then the fact that $\pi_X$ is open follows from Lemma \ref{open} applied to the cartesian left vertical square of the above diagram (taking into account that
$u_X$ is a finite surjective morphism). Thus claim 1 is proved.

\medskip

{\em Claim} 2. The variety $\overline{X}$ is irreducible.

\medskip

Claim 2 follows from Lemma \ref{irrid} applied to the open morphism $\pi_X\colon\overline{X}\to X$, whose fibers are all isomorphic to $\mathbb G_m\times\mathbb G_m$, and in particular, irreducible.

\medskip

By claim 2, passing to the normalization of the irreducible variety $\overline{X}$ we may assume that $\overline{X}$ is normal.
Let $g:Z\to Y$ be the normalization  of $Y=\overline{U}_{X'}$ in the field $K(\overline{X})$ of rational functions of $\overline{X}$ (which makes sense  because the dominant morphism $\overline{X}\to Y$ yields the finite field extension $K(Y)=K(U_{X'})\to K(\overline{X})$). Then we get a commutative diagram of the form   
\begin{diagram}
\overline{X}&\rTo^{i'}&Z\\
\dTo ^{\overline{f}}&  &\dTo _{g}\\
U_{X'}&\rTo^i&Y
\end{diagram}
in which $i$ and $i'$ are open immersions and $g$ is a finite surjective  morphism. Since 
$Y\cap (H_m(e)\times H_n(h))=U_{X'}\cap (H_m(e)\times H_n(h))$, then 
$$\overline{f}^{-1}(U_{X'}\cap (H_m(e)\times H_n(h)))=g^{-1}(U_{X'}\cap (H_m(e)\times H_n(h)))=
g^{-1}(Y\cap (H_m(e)\times H_n(h))).$$
Then by \eqref{final'} and Corollary \ref{universal} we infer that 
$$g^{-1}(Y\cap (H_m(e)\times H_n(h)))=\overline{f}^{-1}(U_{X'}\cap (H_m(e)\times H_n(h)))$$ is connected if $\dim(f(X))>a$ (resp. non-empty if $\dim(f(X))\geq a$). Finally, since $f^{-1}(\Delta)\cong\overline{ f}^{-1}(U_{X'}\cap 
(H_m(e)\times H_n(h)))$, we conclude the proof of Theorem \ref{CT}.\qed

\begin{corollary}\label{CT'}
Let $f:X\to P\times P$ be a morphism as in Theorem $\ref{CT}$ $($i.e. with $X$  a complete irreducible variety and  $P:=\pn m(e)\times \pn n(h)$, $m\geq n\geq 1)$  such that $\dim(f(X))>2m+ n$ $($resp. $\dim(f(X))\geq 2m+n)$. Then 
$f^{-1}(\Delta)$ is connected $($resp. non-empty$)$.\end{corollary}

\proof Since $\dim(X_{13})\leq 2m$, $\dim(X_{24})\leq 2n$ and $m\geq n$, then
$a\leq\max\{2m+n,2n+m\}= 2m+n$, and the conclusion follows from Theorem \ref{CT}. \qed

\medskip
\medskip

We conclude with a few applications of Corollary \ref{CT'} (whose proofs are as in  \cite{FL}).

\begin{corollary}\label{corCT1}
Let $f_i:Y_i\to P=\mathbb P^m(e)\times\mathbb P^n(h)$, $i=1,2$, be two proper morphisms from the complete irreducible varieties $Y_1$ and $Y_2$, with 
 $m\geq n\geq 1$. Assume that $\dim (f_1(Y_1))+\dim (f_2(Y_2))>2m+n$ $($resp. 
$\dim (f_1(Y_1))+\dim (f_2(Y_2)) \geq 2m+n)$. Then $Y_1\times _P Y_2$ is connected $($resp. non-empty$)$.
In particular, if  $Y_1$ and $Y_2$ are closed irreducible subvarieties of $P=\pn m(e)\times \pn n(h)$ with $\dim (Y_1)+\dim (Y_2)>2m+n$ $($resp. $\dim (Y_1)+\dim (Y_2)\geq 2m+n)$, then 
$Y_1\cap Y_2$ is connected $($resp. non-empty$)$.
\end{corollary}

\begin{corollary}\label{corCT2}
Let $Y$ be a projective irreducible variety and let $f:Y\to\pn m(e)\times \pn n(h)$  be a finite unramified morphism, with $m\geq n\geq 1$. If $\dim(Y)>\frac{2m+n}{2}$ then $f$ is a closed embedding.
\end{corollary}

\begin{corollary}\label{corCT3}
Every closed irreducible subvariety $Y$ of $\pn m(e)\times \pn n(h)$ $(m\geq n\geq 1)$   of codimension 
$<\frac{n}{2}$ is algebraically simply connected, i.e. every finite \'etale morphism 
$u:Y'\to Y$,  with $Y'$ connected, is an isomorphism.
\end{corollary}

\itemsep=\smallskipamount
$$\begin{tabular}{llllll}
{\small\it Lucian B\u adescu} & & & & & {\small\it Flavia Repetto}\\
{\small\it Universit\`a degli Studi di Genova}& & & & & {\small\it Universit\`a degli Studi di Milano} \\
{\small\it Dipartimento di Matematica } & & & & & {\small\it Dipartimento di Matematica}\\
{\small\it Via Dodecaneso 35, 16146 Genova, Italy}& & & & & {\small\it Via Saldini 50, 20133 Milano, Italy} \\
{\small\it e-mail: badescu@dima.unige.it} & & & & & {\small\it e-mail: flaviarepetto@yahoo.it}\\
\end{tabular} $$

\end{document}